\documentclass[11pt]{article}
\usepackage{txfonts}
\usepackage{color,epsfig}
\usepackage{mathrsfs}
\usepackage{amssymb}
\usepackage{amsfonts,graphicx,psfrag}

\newtheorem{Th}{Theorem}[section]

\newtheorem{lemma}{Lemma}[section]
\newtheorem{theorem}[lemma]{Theorem}
\newtheorem{corollary}[lemma]{Corollary}

\newtheorem{definition}[lemma]{Definition}

\newtheorem{remark}[lemma]{Remark}
\let\lutzremark=\remark
\def\remark{\lutzremark\normalfont}

\newtheorem{notation}[lemma]{Notation}

\def\be{\begin{equation}}
\def\ee{\end{equation}}
\def\bea{\begin{eqnarray}}
\def\eea{\end{eqnarray}}
\def\bes{\begin{eqnarray*}}
\def\ees{\end{eqnarray*}}

\def\nn{\nonumber}
\def\<{\langle}
\def\>{\rangle}
\def\lb{\label}
\def\bs{\setminus}
\def\pt{\partial}

\def\R{{\bf R}}
\def\C{{\bf C}}
\def\Z{{\bf Z}}
\def\N{{\bf N}}
\def\U{{\bf U}}

\def\F{{\bf F}}

\def\aa{{\alpha}}
\def\bb{{\beta}}
\def\ga{{\gamma}}
\def\Ga{{\Gamma}}

\def\th{{\theta}}
\def\Th{{\Theta}}
\def\om{{\omega}}
\def\Om{{\Omega}}

\def\lm{{\lambda}}

\def\sg{{\sigma}}

\def\A{{\cal A}}

\def\P{{\cal P}}

\def\diag{{\rm diag}}

\def\span{{\rm span}}
\def\Sp{{\rm Sp}}

\def\dm{{\rm \diamond}}

\def\hb{\vrule height0.18cm width0.14cm $\,$}
\def\ol#1{\overline{#1}}
\def\td#1{\tilde{#1}}

\title{The analytical aspect to the linear stability of elliptic equilibrium points
of the Robe's restricted three-body problem}

\author{Qinglong Zhou$^{1} $\thanks{Partially supported by NSFC (No.11501330, No.11425105) of China
           and China Postdoctoral Science Foundation (Grant No. 2015M582071).
           E-mail:zhouqinglong@sdu.edu.cn}\quad
        Yongchao Zhang$^{2} $\thanks{Partially supported by by the Fundamental Research Funds for the Central Universities (Grant No.~N142303010). E-mail: ldfwq@163.com}, \\ \\
$^{1}$ School of Mathematics\\Shandong University, Jinan 250100, Shandong, China\\
$^{2}$ School of Mathematics and Statistics\\Northeastern University at Qinhuangdao 066004, Hebei, China\\}
\date{}

\begin{document}

\maketitle

\begin{abstract}
{We study the Robe's restricted three-body problem.
Such a motion was firstly studied by A. G. Robe in \cite{Robe}, which is used to model
small oscillations of the earth's inner core taking into account the moon attraction.
For the linear stability of elliptic equilibrium points
of the Robe's restricted three-body problem,
earlier results of such linear stability problem depend on a lot of numerical computations,
while we give an analytic approach to it.
The linearized Hamiltonian system near the elliptic relative equilibrium point in our problem
coincides with the linearized system near the Euler elliptic relative equilibria
in the classical three-body problem except for the rang of the mass parameter.
We first establish some relations from the linear stability problem to
symplectic paths and its corresponding linear operators.
Then using the Maslov-type $\omega$-index
theory of symplectic paths and the theory of linear operators, we compute $\om$-indices and obtain
certain properties of the linear stability of elliptic equilibrium points of the Robe's restricted three-body
problem.}
\end{abstract}



\renewcommand{\theequation}{\thesection.\arabic{equation}}

\setcounter{equation}{0}
\setcounter{figure}{0}
\section{Introduction and main results}
\label{sec:1}

A new kind of restricted three-body problem that incorporates the effect of buoyancy
forces was introduced by Robe in 1977.
In \cite{Robe}, he regarded
one of the primaries as a rigid spherical shell $m_1$
filled with a homogenous incompressible fluid of density $\rho_1$.
The second primary is a mass point $m_2$ outside the shell
and the third body $m_3$ is a small solid sphere of density $\rho_3$, inside the shell,
with the assumption that the mass and radius of $m_3$ are infinitesimal.
He has shown the existence of an equilibrium point with $m_3$ at the center of the shell,
where $m_2$ describes a Keplerian orbit around it, see Figure 1.

\begin{figure}[ht]
\centering
\includegraphics[height=4.5cm]{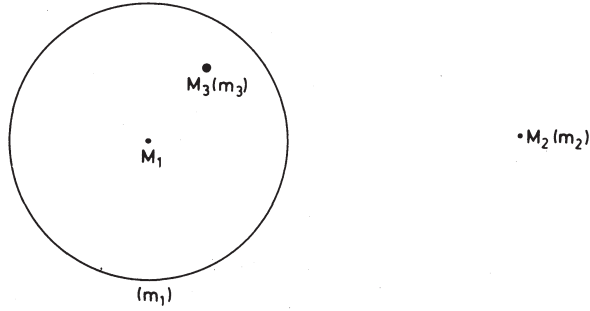}
\caption{The Robe's restricted three-body problem considered:
$m_1$ is a spherical shell filled with a fluid of density $\rho_1$;
$m_2$ a mass point outside the shell and $m_3$ a small solid sphere of density $\rho_3$ inside the shell.}
\end{figure}
\vspace{2mm}

Further, he discussed two cases of the linear stability of the equilibrium points of such restricted three-body problem.
In the first case, the orbit of $m_2$ around $m_1$ is circular and in the second case,
the orbit is elliptic, but the shell is empty (that is no fluid inside it)
or the densities of $m_1$ and $m_3$ are equal.
In the second case, we use ``{\it elliptic equilibrium point}" to call the equilibrium point.
In each case, the domain of stability has been investigated for the whole range
of parameters occurring in the problem.

Later on, A. R. Plastino and A. Plastino (\cite{PlP}) studied the linear stability of the equilibrium points and
the connection between the effect of the buoyancy forces and a perturbation of a Coriolis force.
In 2001, P. P. Hallen and N. Rana (\cite{HaR}) found other new equilibrium points of the restricted problem
and discussed their linear stabilities.
K. T. Singh, B. S. Kushvah and B. Ishwar (\cite{SKI}) examined the stability of triangular equilibrium points in Robe's
generalized restricted three body problem where the problem is generalized in the sense that a more
massive primary has been taken as an oblate spheroid.

However, in \cite{Robe}, for the elliptic equilibrium points,
the bifurcation diagram of linear stability was obtained by numerical methods.
In \cite{PlP,HaR,ShG,SKI}, the authors studied the linear stability of some kinds of equilibrium points,
but their studies did not contain the elliptic case.

On the other hand, in \cite{HS1,HS2} of 2009--2010, X.~Hu and S.~Sun found
a new way to relate the stability problem to the
iterated Morse indices. Recently, by observing new phenomenons and discovering new properties of elliptic Lagrangian
solutions, in the joint paper \cite{HLS} of X. Hu, Y. Long and S. Sun, the linear stability of elliptic Lagrangian solutions
is completely solved analytically by index theory (cf. \cite{Lon1}) and the new results are related directly to
$(\beta,e)$ in the full parameter rectangle.
Inspired by the analytic method, Q. Zhou and Y. Long in \cite{ZL0} studied
the linear stability of elliptic triangle solutions of the charged three-body problem.

Recently, in \cite{ZL,ZL2},
Q. Zhou and Y. Long studied the linear stability of elliptic Euler-Moulton solutions of $n$-body problem
for $n=3$ and for general $n\ge4$, respectively.
Also, the linear stability of Euler collision solutions of $3$-body problem was studied by X. Hu and
Y. Ou in \cite{HO}.

In the current paper,
we study an analytical approach
to the linear stability of equilibrium points of the Robe's restricted three-body problem.
We related their linear stabilities to
the Maslov-type and Morse indices of them.
For such elliptic equilibrium points,
we use index theory to compute the Maslov-type indices of the corresponding symplectic paths
and determine their stability properties.

Following Robe's notations in \cite{Robe}, various forces acting on $m_3$, these are:

(1) The attraction of $m_2$,

(2) The gravitational force of attraction
\begin{equation}
\F_A=-(\frac{4}{3}\pi)G\rho_1m_3\overline{M_1M_3};
\end{equation}
exerted by the fluid of density $\rho_1$;
where $M_1$ is the center of the spherical shell $m_1$, $M_3$ is the center of $M_3$,
and $\overline{M_1M_3}$ is the length of the segment between $M_1$ and $M_3$;

(3) The buoyancy force
\begin{equation}
\F_B=(\frac{4}{3}\pi)G\rho_1^2m_3\overline{M_1M_3}/\rho_3
\end{equation}
of fluid density $\rho_1$.

Let the orbital plane of $m_2$ around $m_1^*$ (that is the shell with its fluid)  be taken as the $x-y$ plane
and let the origin of the coordinate system be at the center of the mass, O, of the two primaries.
The equation of motion of $m_3$ is
\begin{equation}
\ddot{\R}_3=Gm_2\frac{\R_{32}}{R_{32}^2}-\frac{4}{3}\pi G\rho_1(1-\frac{\rho_1}{\rho_3})\R_{13}
\end{equation}
where $\R_3=\overline{OM_3}$ and $\R_{ij}=\overline{M_iM_j}$.

After a detailed calculations,
Robe obtained the equations of the motion:
\begin{eqnarray}
\ddot{x}-2\dot{y}&=&(1+e\cos\theta)^{-1}V_x,
\\
\ddot{y}+2\dot{x}&=&(1+e\cos\theta)^{-1}V_y,
\\
\ddot{z}+\dot{z}&=&(1+e\cos\theta)^{-1}V_z,
\end{eqnarray}
where $\theta$ is the true anomaly in the two-body problem $m_1^*$ and $m_2$,
and $V$ is given by
\begin{equation}
V=\frac{1}{2}(x^2+y^2+z^2)+\frac{\mu}{(x-x_2)+y^2+z^2}
  -\frac{K}{2}(\frac{1-e^2}{1+e\cos\theta})^3[(x-x_1)+y^2+z^2]
\end{equation}
with
\begin{equation}
\mu=\frac{m_2}{m_1^*+m_2},0<\mu<1;\quad K=\frac{4}{3}\pi\frac{\rho_1a^3}{m_1^*+m_2}(1-\frac{\rho_1}{\rho_3}),
\end{equation}
$x_1$ and $x_2$ being the $x$ coordinates of $M_1$ and $M_2$:
\begin{equation}
x_1=\mu,\quad\quad x_2=1-\mu.
\end{equation}

In \cite{Robe}, H. Robe firstly studied the equilibrium points of the problem.
He obtained two kind of equilibrium points,
one is the circular case, and the other is the elliptic case under $K=0$.
He also studied the linear stability of the above two kinds of equilibrium points.
But for the elliptic equilibrium points, only numerical results was obtained.
Later on, in \cite{HaR}, P. P. Hallen and N. Rana studied the existence of all the equilibrium points
in the Robe's restricted three-body problem.
They found that, in the case of equilibrium points with circular, there  are serval
different situations depending on $K$, and the linear stability of such equilibrium points was carefully studied.
More details can be seen in \cite{HaR}.

We focus on the elliptic case when there is no fluid inside the shell or when $\rho_1=\rho_3$,
i.e., $K=0$.
By (17)-(19) in \cite{Robe},the linearized the equations of motion around this equilibrium are:
\begin{eqnarray}
\ddot{x}-2\dot{y}&=&\left\{\frac{1+2\mu}{1+e\cos\theta}-\frac{K(1-e^2)^3}{(1+e\cos\theta)^4}\right\}x,\lb{l1}
\\
\ddot{y}+2\dot{x}&=&\left\{\frac{1-\mu}{1+e\cos\theta}-\frac{K(1-e^2)^3}{(1+e\cos\theta)^4}\right\}y,\lb{l2}
\\
\ddot{z}+\dot{z}&=&\left\{\frac{1-\mu}{1+e\cos\theta}-\frac{K(1-e^2)^3}{(1+e\cos\theta)^4}\right\}z,\lb{l3}
\end{eqnarray}
which is a set of linear homogeneous equations with periodic coefficients of periodic $2\pi$.

Now we study equations (\ref{l1})-(\ref{l3}) from another point of view.
We mainly focus on the linear stability problem on the horizontal plane, i.e., the $xy$-plane,
so we just consider the first two equations.
The linear stability problem along the $z$-axis will be studied in another paper.
Let $(W_1,W_2,w_1,w_2)^T=(\dot{x}-y,\dot{y}+x,x,y)^T$ and $t=\th,K=0$, then we have
\begin{equation}
\frac{d}{dt}
 \left(
 \matrix{
  W_1\cr
  W_2\cr
  w_1\cr
  w_2\cr
 }
 \right)
=
 \left(
 \matrix{
  0 & 1 & -1+\frac{1+2\mu}{1+e\cos{t}} & 0\cr
  -1 & 0 & 0 & -1+\frac{1-\mu}{1+e\cos t}\cr
  1 & 0 & 0 & 1\cr
  0 & 1 & -1 & 0\cr
 }
 \right)
 \left(
 \matrix{
  W_1\cr
  W_2\cr
  w_1\cr
  w_2\cr
 }
 \right).
 \lb{system1}
\end{equation}
Let
\begin{equation}
B(t)=
 \left(
 \matrix{
  1 & 0 & 0 & 1\cr
  0 & 1 & -1 & 0\cr
  0 & -1 & 1-\frac{1+2\mu}{1+e\cos t} & 0\cr
  1 & 0 & 0 & 1-\frac{1-\mu}{1+e\cos t}\cr
 }
 \right),
\label{Bt}
\end{equation}
then (\ref{system1}) can be written as
\begin{equation}
\dot{w}=JB(t)w, \lb{system2}
\end{equation}
where $w=(W_1,W_2,w_1,w_2)^T$.
When $\mu=\bb+1$, $B(t)$ of (\ref{Bt}) coincides with $B(t)$ of (2.35) in \cite{ZL}.
Thus a lot of results which developed in \cite{ZL} can be applied to this paper.

Let $\ga_{\mu,e}(t)$ be the fundamental solution of the linearized Hamiltonian system (\ref{system2}).
Denote by $\Sp(2n)$ the symplectic group of real $2n\times 2n$ matrices. For any
$\omega\in\U=\{z\in\C\;|\;|z|=1\}$ and $M\in\Sp(2n)$, let $\nu_{\om}(M)=\dim_{\C}\ker_{\C}(M-\om I_{2n})$,
and $M$ is called $\omega$-{\it degenerate} ($\om$-{\it non-degenerate} respectively) if $\nu_{\om}(M)> 0$
($\nu_{\om}(M)=0$ respectively). When $\om=1$ and if there is no confusion, we shall simply omit the
subindex $1$ and say just {\it degenerate} or {\it non-degenerate}.

The following two theorems describe main results proved in this paper.

\begin{theorem}\label{T1.1}
In the Robe's restricted three-body problem,
we denote by $\ga_{\mu,e}:[0,2\pi] \to \Sp(4)$ the fundamental solution
of the linearized Hamiltonian system near the equilibrium point.
Then $\ga_{\mu,e}(2\pi)$ is non-degenerate for all $(\mu,e)\in(0,1)\times[0,1)$;
and when $\mu=0$ or $1$, it is degenerate.
Note that the Maslv-type index satisfies $i_1(\ga_{\mu,e})=0$ for all $(\mu,e)\in\Th=[0,1]\times[0,1)$.
Moreover, the following results on the Maslov-type indices of $\ga_{\bb,e}$ hold.

(i) For $e\in[0,1)$, we have
\begin{eqnarray}
&& i_1(\ga_{\mu,e}) = 0. \lb{1.6}\\
&& \nu_1(\ga_{\mu,e}) = \left\{\matrix{
   2, &  {\it if}\;\;\mu=0, \cr
   0, &  {\it if}\;\;0<\mu<1, \cr
   3, &  {\it if}\;\;\mu=1. \cr}\right. \lb{1.7}
\end{eqnarray}

(ii) Let
\be
\mu^*=\frac{5+\sqrt{97}}{16},
\ee
we have
\begin{eqnarray}
&& i_{-1}(\ga_{\mu,0}) = \left\{\matrix{
   0, &  {\it if}\;\;\mu\in[0,\mu^*], \cr
   2, &  {\it if}\;\;\mu\in(\mu^*,1], \cr}\right. \lb{1.9}\\
&& \nu_{-1}(\ga_{\mu,0}) = \left\{\matrix{
   2, &  {\rm if}\;\;\mu=\mu^*, \cr
   0, &  {\rm if}\;\;\mu\in[0,1]\backslash\{\mu^*\}.\cr} \right.
                                   \lb{1.10}
\end{eqnarray}

(iii) For fixed $e\in [0,1)$ and $\om\in\U\backslash\{1\}$, $i_{\om}(\ga_{\mu,e})$ is non-decreasing
and tends from $0$ to $2$
when $\mu$ increases from $0$ to $1$.
\end{theorem}

\begin{remark}\label{R1.2}
(i) Here we are specially interested in indices in eigenvalues $1$ and $-1$. The reason is that the
major changes of the linear stability of the elliptic Euler solutions happen near the eigenvalues
$1$ and $-1$, and such information is used in the next theorem to get the separation curves of the
linear stability domain $[0,1]\times [0,1)$ of the mass and eccentricity parameter $(\mu,e)$.
\end{remark}
\begin{theorem}\label{T1.3}
Using notations in Theorem \ref{T1.1},
for every $e\in[0,1)$, the $-1$-index $i_{-1}(\ga_{\mu,e})$ is non-decreasing,
and strictly decreasing only on two values of $\mu=\mu_1(e)$ and $\mu=\mu_2(e)\in(0,1)$.
Define $\Ga_i=\{(\mu_i(e),e)|e\in[0,1)\}$ for $i=1,2$ and
\be
\mu_m(e)=\min\{\mu_1(e,-1),\mu_2(e,-1)\},\quad
\mu_r(e)=\max\{\mu_1(e,-1),\mu_2(e,-1)\},
\quad e\in[0,1).
\ee

For every $e\in[0,1)$, we also define
\be\lb{mu_l}
\mu_l(e)=\inf\{\mu'\in[0,1]|\sg(\ga_{\mu,e}(2\pi))\cap\U\ne\emptyset,\;\forall \mu\in(0,\mu']\},
\ee
and
\be
\Ga_l=\{(\mu_l(e),e)\in[0,1]\times[0,1)|e\in[0,1)\}.
\ee
Then $\Ga_l$, $\Ga_m$ and $\Ga_r$ from three curves which possess the following properties.

(i) $0<\mu_i(e)<1,i=1,2$ and both $\mu=\mu_1(e)$ and $\mu=\mu_2(e)$ are real analytic in $e\in[0,1)$.
Moreover, $\mu_1(0)=\mu_2(0)=\mu^*$ and $\lim_{e\rightarrow1}\mu_1(e)=\lim_{e\rightarrow1}\mu_2(e)=1$;

(ii)
The two curves $\Ga_1$ and $\Ga_2$ are real analytic in $e$, and bifurcation out  from $(\mu^*,0)$ with tangents
$-\frac{291+15\sqrt{97}}{3104}$ and $\frac{291+15\sqrt{97}}{3104}$ respectively,
thus they are different and their intersection points must be isolated if there exist when $e\in(0,1)$.
Consequently, $\Ga_m$ and $\Ga_r$ are different piecewise real analytic curves, see Figure 2;

(iii) We have
\be
i_{-1}(\ga_{\mu,e}) = \left\{\matrix{
   0, &  {\it if}\;\;\mu\in[0,\mu_m(e)], \cr
   1, &  {\it if}\;\;\mu\in(\mu_m(e),\mu_r(e)],\cr
   2, &  {\it if}\;\;\mu\in(\mu_r(e),1], \cr}\right.
\ee
and $\Ga_m$ and $\Ga_r$ are precisely the $-1$-degenerate curves of the matrix $\ga_{\mu,e}(2\pi)$
in the $(\mu,e)$ rectangle $[0,1]\times[0,1)$;

(iv) Every matrix $\ga_{\mu,e}(2\pi)$ is hyperbolic when $\mu\in[0,\mu_k(e))$, $e\in[0,1)$, and there holds
\be
\mu_l(e)=\sup\{\mu\in[0,1]|\sg(\ga_{\mu,e}(2\pi))\cap\U=\emptyset,\;\forall e\in[0,1)\},
\ee
Consequently, $\Ga_l$ is the boundary curve of the hyperbolic region of $\ga_{\mu,e}(2\pi)$
in the $(\mu,e)$ rectangle $[0,1]\times[0,1)$;

(v) $\Ga_l$ is continuous in $e\in[0,1)$, and $\lim_{e\rightarrow1}\mu_l(e)=1$;

(vi) $\Ga_l$ is different from the curve $\Ga_m$ at least when $e\in[0,\tilde{e})$ for some $\tilde{e}\in(0,1)$;

(vii) We have $\ga_{\mu,e}(2\pi)\approx R(\th_1)\diamond R(\th_2)$ for some $\th_1\in(\pi,2\pi)$ and $\th_2\in(0,\pi)$,
and thus it is strongly linear stable on the segment $\mu_l(e)<\mu<\mu_m(e)$;

(viii) We have $\ga_{\mu,e}(2\pi)\approx R(\th)\diamond D(-2)$ for some $\th\in(\pi,2\pi)$,
and thus it is linearly unstable on the segment $\mu_m(e)<\mu<\mu_r(e)$;

(ix) We have $\ga_{\mu,e}(2\pi)\approx R(\th_1)\diamond R(\th_2)$ for some $\th_1,\th_2\in(\pi,2\pi)$,
and thus it is strongly linear stable on the segment $\mu_r(e)<\mu<1$.
\end{theorem}

\begin{figure}[ht]
\centering
\includegraphics[height=8.5cm]{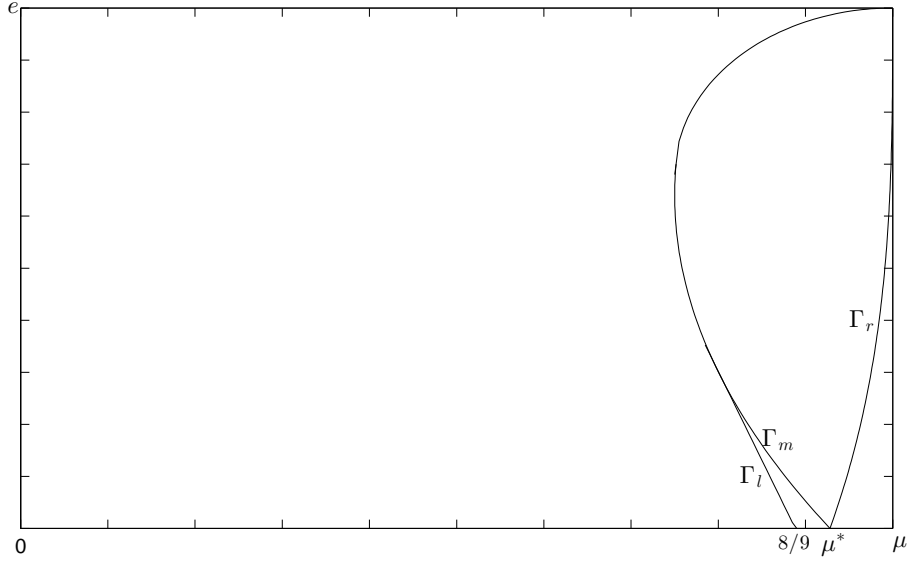}
\caption{Stability bifurcation diagram of elliptic equilibrium points of the Robe's restricted three-body problem
in the $(\mu,e)$ rectangle $[0,1]\times [0,1)$}
\end{figure}
\vspace{2mm}

\begin{remark}\label{R1.4}
For $(\mu,e)$ located on these three special curves, we have the following:

\indent(i) If $\mu_l(e)<\mu_m(e)\le\mu_r(e)$,
we have $\ga_{\mu_l(e),e}(2\pi)\approx N_2(e^{\sqrt{-1}\th},b)$
for some $\th\in(0,\pi)$ and $b=\left(\matrix{b_1&b2\cr b_3&b_4}\right)$ satisfying $(b_2-b_3)\sin\th>0$.
Consequently, the matrix $\ga_{\mu_l(e),e}(2\pi)$ is spectrally stable and linear unstable;

\indent(ii) If $\mu_l(e)=\mu_m(e)<\mu_r(e)$,
we have $\ga_{\mu_l(e),e}(2\pi)\approx N_1(-1,1)\dm D(-2)$ and it is linearly unstable,
or $\ga_{\mu_l(e),e}(2\pi)\approx M_2(-1,c)$ with $c_1,c_2\in\R,c_2\ne0$,
and it is spectrally stable and linearly unstable;

\indent(iii) If $\mu_l(e)=\mu_m(e)=\mu_r(e)$,
we have $\ga_{\mu_l(e),e}(2\pi)\approx N_1(-1,1)\dm N_1(-1,1)$ and it is spectrally stable and linearly unstable;

\indent(iv) If $\mu_l(e)<\mu_m(e)<\mu_r(e)$, we have $\ga_{\mu_m(e),e}(2\pi)\approx N_1(-1,-1)\dm R(\th)$
for some $\th\in(\pi,2\pi)$, and thus is spectrally stable and linearly unstable;

\indent(v) If $\mu_l(e)<\mu_m(e)=\mu_r(e)$, we have $\ga_{\mu_m(e),e}(2\pi)\approx-I_2\dm R(\th)$
for some $\th\in(\pi,2\pi)$, and thus is linearly stable but not strongly linearly stable;

\indent(vi) If $\mu_m(e)<\mu_r(e)$, we have $\ga_{\mu_r(e),e}(2\pi)\approx N_1(-1,1)\dm R(\th)$
for some $\th\in(\pi,2\pi)$, and thus is spectrally stable and linearly unstable.
\end{remark}

The paper is organized as follows. In Section 2, we associate $\ga_{\mu,e}(t)$,
the fundamental solution of the system (\ref{system2}),
with a corresponding second order self-adjoint operator $A(\mu,e)$.
Some connections between $\ga_{\mu,e}(t)$ and $A(\mu,e)$ are given there.
In Section 3, we compute the $\om$-indices along the three boundary segments of $(\mu,e)$ rectangle
$[0,1]\times[0,1)$.
In Section 4, the non-decreasing property of $\om$-index is proved in Lemma \ref{Lemma:increasing.of.index} and Corollary \ref{C4.5}.
Also Theorems 1.1 and Theorem 1.3 are proved there.
For Theorem 1.1, the index properties in
(i)-(iii) are established in Section 3; the non-decreasing property (iv) is proved in Theorem 4.3.
For Theorem 1.3, (i)-(iii) are proved in Section 4.3, and the remain part is proved in Section 4.4.

\section{Associate $\ga_{\mu,e}(t)$ with a second order self-adjoint operator $A(\mu,e)$}

In the Appendix, we give a brief review on the Maslov-type $\om$-index theory for
$\om$ in the unit circle of the complex plane following \cite{Lon4}. In the following, we use notations
introduced there.

Let
\be  J_2=\left(\matrix{ 0 & -1 \cr 1 & 0 \cr}\right), \qquad
    Q_{\mu,e}(t)=\left(\matrix{\frac{1+2\mu}{1+e\cos t} & 0 \cr
                    0 & \frac{1-\mu}{1+e\cos t} \cr}\right),
                                     \lb{2.20}\ee
and set
\be L(t,x,\dot{x})=\frac{1}{2}\|\dot{x}\|^2 + J_2x(t)\cdot\dot{x}(t) + \frac{1}{2}Q_{\mu,e}(t)x(t)\cdot x(t),
       \qquad\quad  \forall\;x\in W^{1,2}(\R/2\pi\Z,\R^2),  \lb{Lag1}\ee
where $a\cdot b$ denotes the inner product in $\R^2$.
By Legendrian transformation, the corresponding
Hamiltonian function to system (\ref{system2}) is
$$   H(t,w)=\frac{1}{2}B(t)w\cdot w,\qquad \forall\; w\in\R^4.  $$

Now let $\ga=\ga_{\mu,e}(t)$ be the fundamental solution of the (\ref{system2}) satisfies:
\begin{eqnarray}
\dot{\ga}(t)&=&JB(t)\ga(t), \lb{path.eq}
\\
\ga(0)&=&I_4.   \lb{path.bc}
\end{eqnarray}
In order to transform the Lagrangian system (\ref{Lag1}) to a simpler linear operator corresponding to
a second order Hamiltonian system with the same linear stability as $\ga_{\mu,e}(2\pi)$, using $R(t)$
and $R_4(t)=\diag(R(t),R(t))$ as in Section 2.4 of \cite{HLS}, we let
\be  \xi_{\mu,e}(t) = R_4(t)\ga_{\mu,e}(t), \qquad \forall\; t\in [0,2\pi],
(\mu,e)\in [0,1]\times [0,1). \lb{2.25}\ee
One can show by direct computation that
\be  \frac{d}{dt}\xi_{\mu,e}(t)
  = J \left(\matrix{I_2 & 0 \cr
                    0 & R(t)(I_2-Q_{\mu,e}(t))R(t)^T \cr}\right)\xi_{\mu,e}(t). \lb{2.26}\ee
Note that $R_4(0)=R_4(2\pi)=I_4$, so $\ga_{\mu,e}(2\pi)=\xi_{\mu,e}(2\pi)$ holds and the linear stabilities
of the systems (\ref{path.bc}) and (\ref{2.26}) are precisely the same.

By (\ref{2.25}) the symplectic paths $\ga_{\mu,e}$ and $\xi_{\mu,e}$ are homotopic to each other via the
homotopy $h(s,t)=R_4(st)\ga_{\mu,e}(t)$ for $(s,t)\in [0,1]\times [0,2\pi]$. Because $R_4(s)\ga_{\mu,e}(2\pi)$
for $s\in [0,1]$ is a loop in $\Sp(4)$ which is homotopic to the constant loop $\ga_{\mu,e}(2\pi)$, we have
$\ga_{\mu,e} \sim_1 \xi_{\mu,e}$ by the homotopy $h$. Then by Lemma 5.2.2 on p.117 of \cite{Lon4}, the
homotopy between $\ga_{\mu,e}$ and $\xi_{\mu,e}$ can be realized by a homotopy which fixes the end point
$\ga_{\mu,e}(2\pi)$ all the time. Therefore by the homotopy invariance of the Maslov-type index (cf. (i) of
Theorem 6.2.7 on p.147 of \cite{Lon4}) we obtain
\be  i_{\om}(\xi_{\mu,e}) = i_{\om}(\ga_{\mu,e}), \quad \nu_{\om}(\xi_{\mu,e}) = \nu_{\om}(\ga_{\mu,e}),
       \qquad \forall \,\omega\in\U, \; (\mu,e)\in [0,1]\times [0,1). \lb{2.27}\ee
Note that the first order linear Hamiltonian system (\ref{2.26}) corresponds to the following second order
linear Hamiltonian system
\be  \ddot{x}(t)=-x(t)+R(t)Q_{\mu,e}(t)R(t)^Tx(t). \lb{2.28}\ee

For $(\mu,e)\in [0,1]\times [0,1)$, the second order differential operator corresponding to (\ref{2.28}) is given by
\bea  A(\mu,e)
&=& -\frac{d^2}{dt^2}I_2-I_2+R(t)Q_{\mu,e}(t)R(t)^T  \nn\\
&=& -\frac{d^2}{dt^2}I_2-I_2+\frac{1}{2(1+e\cos t)}[(2+\mu)I_2+3{\mu}S(t)],
\lb{2.29}
\label{A_mu.e}
\eea
where $S(t)=\left(\matrix{ \cos 2t & \sin 2t \cr
                           \sin 2t & -\cos 2t \cr}\right)$, defined on the domain $\ol{D}(\omega,2\pi)$
in (\ref{2.11}). Then it is self-adjoint and depends on the parameters $\mu,K$ and $e$. By Lemma
\ref{L2.3}, we have for any $\mu,K$ and $e$, the Morse index $\phi_{\om}(A(\mu,K,e))$ and nullity $\nu_{\om}(A(\mu,K,e))$
of the operator $A(\mu,e)$ on the domain $\ol{D}(\omega,2\pi)$ satisfy
\be  \phi_{\om}(A(\mu,e)) = i_{\om}(\xi_{\mu,e}), \quad \nu_{\om}(A(\mu,e)) = \nu_{\om}(\xi_{\mu,e}), \qquad
           \forall \,\om\in\U. \lb{2.30}\ee

In the rest of this paper, we shall use both of the paths $\ga_{\mu,e}$ and $\xi_{\mu,e}$ to study
the linear stability of $\ga_{\mu,e}(2\pi)=\xi_{\mu,e}(2\pi)$. Because of (\ref{2.27}), in many cases and
proofs below, we shall not distinguish these two paths.

\setcounter{equation}{0}
\section{Computation of the $\om$-indices on the boundary of the bounded rectangle $[0,1]\times [0,1)$}
\label{sec:3}

We first know the full range of $(\mu,e)$ is $(0,1)\times[0,1)$.
For convenience in the mathematical study, we extend the range of $(\mu,e)$ to $[0,1]\times [0,1)$.

Furthermore, we need more precise information on indices and stabilities of $\ga_{\mu,e}$ at the boundary of the
$(\mu,e)$ rectangle $[0,1]\times [0,1)$.

\subsection{$\om$-indices on the boundary segments $\{0\}\times [0,1)$ and $\{1\}\times [0,1)$}

When $\mu=0$, from (\ref{A_mu.e}), we have
\be
A(0,e)=-\frac{d^2}{dt^2}I_2-I_2+\frac{1}{1+e\cos t}I_2,
\ee
this is just the same case which has been discussed in Section 4.1 of \cite{ZL}.
Using Lemma 4.1 of \cite{ZL}, $A(0,e)$ is non-negative definite for the $\om=1$ boundary condition,
and  $A(0,e)$ is positive definite for the $\om\in\U\backslash{1}$ boundary condition.
Hence we have
\begin{eqnarray}
i_{\om}(\ga_{0,e})=i_{\om}(\xi_{0,e})=
\left\{
\begin{array}{l}
0,\quad {\rm if}\;\; \om=1,\\
0,\quad {\rm if}\;\; \om \in \U\bs\{1\},
\end{array}
\right.
\lb{index.of.0e}
\\
\nu_{\om}(\ga_{0,e}) = \nu_{\om}(\xi_{0,e})
= \left\{
\begin{array}{l}
2, \quad {\rm if}\;\;\om = 1,
\\
0, \quad {\rm if}\;\;\om \in \U\bs\{1\}.
\end{array}
\right.
\lb{null.index.of.0e}
\end{eqnarray}

When $\mu=1$, from (\ref{A_mu.e}), we have
\be
A(0,e)=-\frac{d^2}{dt^2}I_2-I_2+\frac{3}{2(1+e\cos t)}(I_2+3S(t)).
\ee
This is just the case which has been discussed in Section 3.1 of \cite{HLS}.
We just cite the results here:
\begin{eqnarray}
i_{\om}(\ga_{1,e})=i_{\om}(\xi_{1,e})=
\left\{
\begin{array}{l}
0,\quad {\rm if}\;\; \om=1,\\
2,\quad {\rm if}\;\; \om \in \U\bs\{1\},
\end{array}
\right.
\lb{index.of.1e}
\\
\nu_{\om}(\ga_{1,e}) = \nu_{\om}(\xi_{1,e})
= \left\{
\begin{array}{l}
3, \quad {\rm if}\;\;\om = 1,
\\
0, \quad {\rm if}\;\;\om \in \U\bs\{1\}.
\end{array}
\right.
\lb{null.index.of.1e}
\end{eqnarray}

\subsection{$\om$-indices on the boundary $[0,1]\times \{0\}$}

In this case $e=0$. It is considered in (A) of Subsection 3.1 of \cite{HLS} when $\bb=0$.
Below, we shall first recall the
properties of eigenvalues of $\ga_{\bb,0}(2\pi)$. Then we carry out the computations of normal
forms of $\ga_{\bb,0}(2\pi)$, and $\pm 1$ indices $i_{\pm 1}(\ga_{\bb,0})$ of the path
$\ga_{\bb,0}$ for all $\bb\in [0,\infty)$, which are new.

In this case, the linearized system (\ref{system1}) becomes an ODE
system with constant coefficients:
\be B = B(t) = \left(\matrix{1 & 0 & 0 & 1\cr
                             0 & 1 & -1 & 0 \cr
                             0 & -1 & -2\mu & 0 \cr
                             1 & 0 & 0 & \mu \cr}\right).  \lb{3.10}\ee
The characteristic polynomial $\det(JB-\lambda I)$ of $JB$ is given by
\be \lambda^4 + (2-\mu)\lambda^2+(1-\mu)(1+2\mu) = 0.  \lb{3.11}\ee
Letting $\aa=\lambda^2$, the two roots of the quadratic polynomial $\aa^2 + (2-\mu)\aa +(1-\mu)(1+2\mu)$
are given by $\aa_1=\frac{\mu-2+\sqrt{9\mu^2-8\mu}}{2}$
and $\aa_2=\frac{\mu-2-\sqrt{9\mu^2-8\mu}}{2}$.
Therefore the four roots of the polynomial
(\ref{3.11}) are given by
\begin{eqnarray}
\aa_{1,\pm} &=&\pm\sqrt{\aa_1}= \pm\sqrt{\frac{\mu-2+\sqrt{9\mu^2-8\mu}}{2}},\lb{3.12a}
\\
\aa_{2,\pm} &=&\pm\sqrt{\aa_2}= \pm\sqrt{\frac{\mu-2-\sqrt{9\mu^2-8\mu}}{2}}. \lb{3.12}
\end{eqnarray}

{\bf (A)  Eigenvalues of $\ga_{\mu,0}(2\pi)$ for $\mu\in [0,1]$.}

When $0<\mu<{8\over9}$, from (\ref{3.12a}) and (\ref{3.12}) by direct computation
the four characteristic multipliers of the matrix $\ga_{\mu,0}(2\pi)$ is given by
\begin{eqnarray}
\rho_{1,\pm}(\beta) = e^{2\pi\aa_{1,\pm}}\in\C\backslash(\U\cap\R),
\\
\rho_{2,\pm}(\beta) = e^{2\pi\aa_{2,\pm}}\in\C\backslash(\U\cap\R).
\lb{rho.CS}
\end{eqnarray}

When ${8\over9}\le\mu\le1$,
by (\ref{3.12a}) and (\ref{3.12}), we get four characteristic multipliers of $\ga_{\mu,0}(2\pi)$
\begin{eqnarray}
\rho_{i,\pm}(\mu) = e^{2\pi\aa_{1,\pm}} =  e^{\pm 2\pi\sqrt{-1}\th_i(\mu)},
\quad i=1,2,
\lb{3.13a}
\end{eqnarray}
where
\be
\th_1(\mu) = \sqrt{\frac{2-\mu-\sqrt{9\mu^2-8\mu}}{2}},\quad
\th_2(\mu) = \sqrt{\frac{2-\mu+\sqrt{9\mu^2-8\mu}}{2}}.
\lb{3.14}
\ee

Moreover, when ${8\over9}\le\mu\le1$, we have
\begin{eqnarray}
\frac{d\th_1(\mu)}{d\mu}=-\frac{1+\frac{9\mu-4}{2\sqrt{9\mu^2-8\mu}}}{4\sqrt{\frac{2-\mu-\sqrt{9\mu^2-8\mu}}{2}}}<0,\label{da1.db}\\
\frac{d\th_2(\mu)}{d\mu}=\frac{-1+\frac{9\mu-4}{2\sqrt{9\mu^2-8\mu}}}{4\sqrt{\frac{2-\mu+\sqrt{9\mu^2-8\mu}}{2}}}>0.\label{da2.db}
\end{eqnarray}
Thus $\th_1(\mu)$ and $\th_2(\mu)$ are monotonic with respect to $\mu$ in this case.

From (\ref{3.14}), $\th_1({8\over9})=\th_2({8\over9})={\sqrt{5}\over3}$ and $\th_1(1)=0,\th_2(1)=1$.
Letting $\mu^*$ be the $\mu$ such that $\th_1(\mu^*)={1\over2}$, then we have
\be
\mu^*=\frac{5+\sqrt{97}}{16}.
\ee
It is obvious that $\mu^*>\frac{8}{9}$.

Specially, we obtain the following results:

When $\mu=0$, we have $\sg(\ga_{0,0}(2\pi)) = \{1, 1, 1, 1\}$.

When $0<\mu<\frac{8}{9}$, using notations defined in (\ref{rho.CS}),
the four characteristic multipliers of $\ga_{\mu,0}(2\pi)$ satisfy
$\sg(\ga_{\mu,0}(2\pi))\in\C\backslash(\U\cap\R)$.

When $\mu=\frac{8}{9}$, we have double eigenvalues
$\rho_{1,\pm}=\rho_{2,\pm}=e^{\pm\sqrt{-1}\frac{2\sqrt{5}}{3}\pi}$.

When $\frac{8}{9}<\mu<\frac{5+\sqrt{97}}{16}(=\mu^*)$,
in (\ref{3.14}), the angle $\th_1(\mu)$ decreases strictly from $\frac{\sqrt{5}}{3}$ to $\frac{1}{2}$,
and the $\th_2(\mu)$ increases strictly from $\frac{\sqrt{5}}{3}$ to ${\sqrt{23-\sqrt{97}}\over4}$
as $\mu$ increases from $\frac{8}{9}$ to $\mu^*$.
Thus specially, we obtain $\sg(\ga_{\mu,0}(2\pi))\in\U\backslash\R$.

When $\mu=\frac{5+\sqrt{97}}{16}(=\mu^*)$, we have $\th_1(\mu^*)={1\over2}$ and
$\th_2(\mu^*)={\sqrt{23-\sqrt{97}}\over4}$.
Therefore we obtain $\rho_{1,\pm}(\hat\bb_{\frac{1}{2}})=e^{\pm \sqrt{-1} \pi} = -1$
and $\rho_{2,\pm}=e^{\pm \sqrt{-1} {\sqrt{23-\sqrt{97}}\over2}\pi}\in\U\backslash\R$.

When $\frac{5+\sqrt{97}}{16}(=\mu^*)<\mu<1$,
the angle $\th_1(\mu)$ decreases strictly from $\frac{1}{2}$ to $0$,
and the $\th_2(\mu)$ increases strictly from ${\sqrt{23-\sqrt{97}}\over4}$ to $1$
as $\mu$ increases from $\frac{8}{9}$ to $\mu^*$.
Thus we obtain $\sg(\ga_{\mu,0}(2\pi))\in\U\backslash\R$.

When $\mu=1$, we have $\th_1(1)=0$ and $\th_2(1)=1$, and hence
$\sg(\ga_{0,0}(2\pi)) = \{1, 1, 1, 1\}$.

{\bf (B)  Indices $i_1(\ga_{\mu,0})$ of $\ga_{\mu,0}(2\pi)$ for $\mu\in [0,1]$.}

Define
\begin{equation}\lb{f0}
f_{0,1}=R(t)\left(\matrix{1\cr 0}\right),\quad
f_{0,2}=R(t)\left(\matrix{0\cr 1}\right),
\end{equation}
and
\begin{equation}\lb{fn}
f_{n,1}=R(t)\left(\matrix{\cos nt\cr 0}\right),\quad
f_{n,2}=R(t)\left(\matrix{0\cr \cos nt}\right),\quad
f_{n,3}=R(t)\left(\matrix{\sin nt\cr 0}\right),\quad
f_{n,4}=R(t)\left(\matrix{0\cr \sin nt}\right),
\end{equation}
for $n\in\N$.
Then $f_{0,1},\;f_{0,2}$ and $f_{n,1},\;f_{n,2}\;f_{n,3},\;f_{n,4}\;n\in\N$
form an orthogonal basis of $\overline{D}(1,2\pi)$.
By (\ref{2.29}) and $\frac{dR(t)}{dt}=JR(t)$, computing $A(\bb,0)f_{n,1}$ yields
\begin{eqnarray}
A(\mu,0)f_{n,1}&=&[-\frac{d^2}{dt^2}I_2-I_2+R(t)K_{0,e}(t)R(t)^T]R(t)\left(\matrix{\cos nt\cr 0}\right)
\nonumber\\
&=&R(t)\left(\matrix{(n^2+1+2\mu)\cos nt\cr 2n\sin nt}\right)
\nonumber\\
&=&(n^2+1+2\mu)f_{n,1}+2nf_{n,4}.
\end{eqnarray}
Similarly, we have
\begin{eqnarray}
\left(\matrix{A(\mu,0) & O\cr O & A(\mu,0)}\right)\left(\matrix{f_{0,1}\cr f_{0,2}}\right)&=&
\left(\matrix{1+2\mu & 0\cr 0 & 1-\mu}\right)
\left(\matrix{f_{0,1}\cr f_{0,2}}\right),\lb{A1}
\\
\left(\matrix{A(\mu,0) & O\cr O & A(\mu,0)}\right)\left(\matrix{f_{n,1}\cr f_{n,4}}\right)&=&
\left(\matrix{n^2+1+2\mu & 2n\cr 2n & n^2+1-\mu}\right)
\left(\matrix{f_{n,1}\cr f_{n,4}}\right),\lb{A2}
\\
\left(\matrix{A(\mu,0) & O\cr O & A(\mu,0)}\right)\left(\matrix{f_{n,3}\cr f_{n,2}}\right)&=&
\left(\matrix{n^2+1+2\mu & -2n\cr -2n & n^2+1-\mu}\right)
\left(\matrix{f_{n,3}\cr f_{n,2}}\right),\lb{A3}
\end{eqnarray}
for $n\in\N$.
Denoting
\begin{equation}
B_0=\left(\matrix{1+2\mu & 0\cr 0 & 1-\mu}\right),\quad
B_n=\left(\matrix{n^2+1+2\mu & 2n\cr 2n & n^2+1-\mu}\right),\quad
\tilde{B}_n=\left(\matrix{n^2+1+2\mu & -2n\cr -2n & n^2+1-\mu}\right).
\end{equation}
Denote the characteristic polynomial of $B_n$ and $\tilde{B}_n$ by $p_n(\lambda)$ and $\tilde{p}_n(\lambda)$
respectively, then we have
\begin{equation}
p_n(\lambda)=\tilde{p}_n(\lambda)=\lambda^2-(2n^2+2+\mu)\lambda-[2\mu^2-(n^2+1)\mu-(n^2-1)^2].
\end{equation}

If $\mu=0$, then $B_0>0$ and $p_n(\lambda)=[\lambda-(n+1)^2][\lambda-(n-1)^2]$, and hence
both $B_1$ and $\tilde{B}_1$ have a zero eigenvalue, and all other eigenvalues of $B_n$ and $\tilde{B}_n(n\ge1)$
are positive. Then we have  $i_1(\ga_{0,0})=0$ and $\nu_1(\ga_{0,0})=2$.

If $0<\mu<1$, then $B_0>0$ and both $B_1$ and $\tilde{B}_1$ have two positive eigenvalues.
Moreover, we have $B_n>\left(\matrix{n^2 & 2n\cr 2n & n^2}\right)>0$,
and hence when $n\ge2$, $B_n$ has two positive eigenvalues.
Similarly, when $n\ge2$, $\tilde{B}_n$ has two positive eigenvalues.
Then we have $i_1(\ga_{0,0})=0$ and $\nu_1(\ga_{0,0})=0$.

Therefore, we have
\begin{eqnarray}
&& i_1(\ga_{\mu,0}) = 0 ;\lb{1-index.of.b0}
\\
&& \nu_1(\ga_{\mu,0}) = \left\{\matrix{
                 2, &  {\rm if}\;\;\mu=0, \cr
                 0, &  {\rm if}\;\;0<\mu<1, \cr
                 3, &  {\rm if}\;\;\mu=1. \cr}\right. \lb{null.1-index.of.b0}
\end{eqnarray}

{\bf (C)  Indices $i_{\om}(\ga_{\mu,0}),\;\om\ne1$ for $\mu\in [0,1]$.}

Because $B(t)$ is a constant matrix depending only on $\mu$ when $e=0$,
it is possible to compute the fundamental matrix path $\ga_{\mu,0}(t)$ explicitly.
Using the notations in (\A), we have $v_{-1}(\ga_{\mu^*,0})$
\begin{eqnarray}
\nu_{-1}(\ga_{\mu,0}) = \left\{\matrix{
                 2, &  {\rm if}\;\;\mu=\frac{5+\sqrt{97}}{16}(=\mu^*), \cr
                 0, &  {\rm if}\;\;\mu\ne\frac{5+\sqrt{97}}{16}. \cr}\right. \lb{null.-1.index.of.b0}
\end{eqnarray}
By a similar analysis to (B), we have
\begin{eqnarray}
&& i_{-1}(\ga_{\mu,0}) = \left\{\matrix{  
                 0, &  {\rm if}\;\;0\le\mu\le\mu^*, \cr
                 2, &  {\rm if}\;\;\mu^*<\mu\le1, \cr}\right.
\lb{-1.index.of.b0}
\end{eqnarray}

\setcounter{equation}{0}
\section{The separation curves of the different linear stability patterns of the elliptic equilibrium points through different $(\mu,e)$ parameters}\label{sec:4}

\subsection{The increasing of $\om$-indeces of $\ga_{\mu,e}$}

For $(\mu,e)\in (0,1)\times [0,1)$, we can rewrite $A(\mu,e)$ as follows
\begin{eqnarray} A(\mu,e) &=& -\frac{d^2}{dt^2}I_2-I_2+\frac{I_2}{1+e\cos t}+\frac{\mu}{2(1+e\cos t)}(I_2+3S(t))
\nonumber\\
&=& \mu\bar{A}(\mu,e), \lb{4.11}
\end{eqnarray}
where we define
\be\label{bar.A}
\bar{A}(\mu,e)=\frac{-\frac{d^2}{dt^2}I_2-I_2+\frac{I_2}{1+e\cos t}}{\mu}+\frac{I_2+3S(t)}{2(1+e\cos t)}
=\frac{A(0,e)}{\mu}+\frac{I_2+3S(t)}{2(1+e\cos t)}.
\ee

Note that $A(\mu,e)$ is the same operator of $A(\bb,e)$ when $\bb=\mu-1$ in \cite{ZL}
with different parameter ranges,
by Lemma 4.2 in \cite{ZL} and modifying its proof to the different range of parameters,
we get the following important lemma:

\begin{lemma}\label{Lemma:increasing.of.index}
(i) For each fixed $e\in [0,1)$, the operator $\bar{A}(\mu,e)$ is non-increasing
with respect to $\mu\in (0,1)$ for any fixed $\omega\in\U$. Specially
\be  \frac{\pt}{\pt\mu}\bar{A}(\mu,e)|_{\mu=\mu_0} = -\frac{1}{\mu_0^2}A(0,e),  \lb{4.14}\ee
is a non-negative definite operator for $\mu_0\in (0,1)$.

(ii) For every eigenvalue $\lm_{\mu_0}=0$ of $\bar{A}(\mu_0,e_0)$ with $\om\in\U$ for some
$(\mu_0,e_0)\in (0,1)\times [0,1)$, there holds
\be \frac{d}{d\mu}\lm_{\mu}|_{\mu=\mu_0} < 0.  \lb{4.15}\ee

(iii)  For every $(\mu,e)\in(0,1)\times[0,1)$ and $\om\in\U$,
there exist $\epsilon_0=\epsilon_0(\mu,e)>0$ small enough such that for all $\epsilon\in(0,\epsilon_0)$ there holds
\be
i_\om(\ga_{\mu+\epsilon,e})-i_\om(\ga_{\mu,e})=\nu_\om(\ga_{\mu,e}).
\ee
\end{lemma}

Consequently we arrive at
\begin{corollary}\label{C4.5} For every fixed $e\in [0,1)$ and $\om\in \U$, the index function
$\phi_{\om}(A(\mu,e))$, and consequently $i_{\om}(\ga_{\mu,e})$, is non-decreasing
as $\mu$ increases from $0$ to $1$.
When $\om=1$, these index functions are constantly equal to $0$,
and when $\om\in\U\bs\{1\}$, they are increasing and tends from $0$ to $2$.
\end{corollary}

{\bf Proof.}
For $0<\mu_1<\mu_2\le1$ and fixed $e\in [0,1)$, when $\mu$ increases from $\mu_1$ to
$\mu_2$, it is possible that positive eigenvalues of $\bar{A}(\mu_1,e)$ pass through $0$ and become
negative ones of $\bar{A}(\mu_2,e)$, but it is impossible that negative eigenvalues of
$\bar{A}(\mu_2,e)$ pass through $0$ and become positive by (ii) of Lemma \ref{Lemma:increasing.of.index}.
\hb

\subsection{The $\om$-degenerate curves of of $\ga_{\mu,e}$}

By a similar analysis to the proof of Proposition 6.1 in \cite{HLS},
for every $e\in[0,1)$ and $\om\in\U\backslash\{1\}$,
the total multiplicity of $\om$-degeneracy of $\ga_{\mu,e}(2\pi)$ for $\mu\in[0,1]$ is always precisely 2, i.e.,
\be
\sum_{\mu\in[0,1]}v_\om(\ga_{\mu,e}(2\pi))=2,\quad \forall \om\in\U\backslash\{1\}.
\ee

Consequently,
together with the definiteness of $A(0,e)$ for the $\om\in\U\backslash\{1\}$ boundary condition,
we have
\begin{theorem}\label{Th:om.degenerate.curves}
For any $\om\in\U\backslash\{1\}$, there exist two analytic $\om$-degenerate curves $(\mu_i(e,\om),e)$
in $e\in[0,1)$ with $i=1,2$.
Specially, each $\mu_i(e,\om)$ is areal analytic function in $e\in[0,1)$,
and $0<\mu_i(e,\om)<1$ and $\ga_{\mu_i(e,\om),e}(2\pi)$ is $\om$-degenerate for $\om\in\U\backslash\{1\}$
and $i=1,2$.
\end{theorem}

{\bf Proof.}
We prove first that $i_\om(\ga_{\mu,e})=0$ when $\mu$ is near $0$.
By Lemma 4.1(ii) in \cite{ZL}, $A(0,e)$ is positive definite on $\overline{D}(\om,2\pi)$.
Therefore, there exists an $\epsilon>0$ small enough such that
$A(\mu,e)$ is also positive definite on $\overline{D}(\om,2\pi)$ when $0<\mu<\epsilon$.
Hence $\nu_\om(\ga_{\mu,e})=\nu_\om(A(\mu,e))=0$ when $0<\mu<\epsilon$.
Thus we have proved our claim.

Then under a similar steps to those of Lemma 6.2 and Theorem 6.3 in in \cite{HLS},
we can prove the theorem.
\hb

\subsection{The $\omega=-1$ degenerate curves of $\ga_{\mu,e}$}

Specially, for $\om=-1$, $e\in[0,1)$ we define
\be
\mu_m(e)=\min\{\mu_1(e,-1),\mu_2(e,-1)\},\quad
\mu_r(e)=\max\{\mu_1(e,-1),\mu_2(e,-1)\},
\ee
where $\mu_i(e,-1)$ are the two $-1$-dgenerate curves as in Theorem \ref{Th:om.degenerate.curves}.

By (\ref{null.-1.index.of.b0}), $-1$ is a double eigenvalue of the matrix $\ga_{\mu^*,e}(2\pi)$,
then the two curves bifurcation out from $(\mu^*,0)$ when $e>0$ is small enough.

Recall $A(\mu^*,0)$ is $-1$-degenerate and by (\ref{null.-1.index.of.b0}),
$\dim\ker A(\mu^*,0)=v_{-1}(\ga_{\mu^*,0})=2$.
By the definition of (\ref{2.11}),
we have $R(t)\left(\matrix{\tilde{a}_n\sin (n+\frac{1}{2})t\cr \cos (n+\frac{1}{2})t}\right)\in\ol{D}(-1,2\pi)$
for any constant $\tilde{a}_n$.

Moreover,
$A(\mu,0)R(t)\left(\matrix{\tilde{a}_n\sin (n+\frac{1}{2})t\cr \cos (n+\frac{1}{2})t}\right)=0$ reads
\be
\left\{
\begin{array}{cr}
(n+\frac{1}{2})^2\tilde{a}_n-2(n+\frac{1}{2})+(1+2\mu)\tilde{a}_n&=0,
\\
(n+\frac{1}{2})^2-2(n+\frac{1}{2})\tilde{a}_n+1-\mu&=0.
\end{array}
\right.
\ee
Then $2\mu^2-((n+\frac{1}{2})^2+1)\mu-[(n+\frac{1}{2})^2-1]^2=0$ which holds only when $n=0$
and $\mu=\frac{5+\sqrt{97}}{16}=\mu^*$ again and
\be\lb{tilde.a}
\tilde{a}_0={1\over4}+1-\mu^*=\frac{15-\sqrt{97}}{16}.
\ee
Then we have
$R(t)\left(\matrix{\tilde{a}_0\sin\frac{t}{2}\cr \cos\frac{t}{2}}\right)\in\ker A(\mu^*,0)$.
Similarly
$R(t)\left(\matrix{\tilde{a}_0\cos\frac{t}{2}\cr -\sin\frac{t}{2}}\right)\in\ker A(\mu^*,0)$,
therefore we have
\begin{equation}
\ker A(\mu^*,0)=\span\left\{
R(t)\left(\matrix{\tilde{a}_0\sin\frac{t}{2}\cr \cos\frac{t}{2}}\right),\quad
R(t)\left(\matrix{\tilde{a}_0\cos\frac{t}{2}\cr -\sin\frac{t}{2}}\right)
\right\}.
\lb{ker.A.of.-1}
\end{equation}
Indeed, we have the following theorem:

\begin{theorem}\label{Th:tangent.direction}
The tangent direction of the two curves $\Ga_m$ and $\Ga_r$ bifurcation from $(\mu^*,0)$
when $e>0$ is small are given by
\be
\mu_m'(e)|_{e=0}=-\frac{291+15\sqrt{97}}{3104},\quad
\mu_r'(e)|_{e=0}=\frac{291+15\sqrt{97}}{3104}.
\ee
\end{theorem}

{\bf Proof.}
Now let $(\mu(e),e)$ be one of such curves
(i.e., one of $(\mu_i(-1,e),e),i=1,2$.)
which starts from $\mu^*$ with $e\in[0,\epsilon)$ for some small $\epsilon>0$ and $x_e\in\bar{D}(1,2\pi)$
be the corresponding eigenvector, that is
\be
A(\mu(e),e)x_e=0.
\ee
Without loose of generality, by (\ref{ker.A.of.-1}), we suppose
$$
z = (\tilde{a}_0\sin\frac{t}{2},\cos\frac{t}{2})^T
$$
and
\be
x_0=R(t)z=R(t)(\tilde{a}_0\sin\frac{t}{2},\cos\frac{t}{2})^T.\lb{4.71}
\ee
There holds
\be
\<A(\mu(e),e)x_e,x_e\>=0.\label{Axx-1}
\ee

Differentiating both side of (\ref{Axx-1}) with respect to $e$ yields
$$ \mu'(e)\<\frac{\pt}{\pt \mu}A(\mu(e),e)x_e,x_e\> + (\<\frac{\pt}{\pt e}A(\mu(e),e)x_e,x_e\>
       + 2\<A(\mu(e),e)x_e,x'_e\> = 0,  $$
where $\mu'(e)$ and $x'_e$ denote the derivatives with respect to $e$. Then evaluating both
sides at $e=0$ yields
\be  \mu'(0)\<\frac{\pt}{\pt \mu}A(\mu^*,0)x_0,x_0\>
   + \<\frac{\pt}{\pt e}A(\mu^*,0)x_0,x_0\> = 0.
\lb{4.73}
\ee
Then by the definition (\ref{2.29}) of $A(\mu,e)$ we have
\bea
\left.\frac{\pt}{\pt\mu}A(\mu,e)\right|_{(\mu,e)=(\mu^*,0)}
    &=& \left.R(t)\frac{\pt}{\pt\mu}K_{\mu,e}(t)\right|_{(\mu,e)=(\mu^*,0)}R(t)^T,  \lb{4.74}\\
\left.\frac{\pt}{\pt e}A(\mu,e)\right|_{(\mu,e)=(\mu^*,0)}
    &=& \left.R(t)\frac{\pt}{\pt e}K_{\mu,e}(t)\right|_{(\mu,e)=(\mu^*,0)}R(t)^T,  \lb{4.75}\eea
where $R(t)$ is given in \S 2.1. By direct computations from the definition of $K_{\mu,e}(t)$ in
(\ref{2.20}), we obtain
\bea
&& \frac{\pt}{\pt\mu}K_{\mu,e}(t)\left|_{(\mu,e)=(\mu^*,0)}
       = \left(\matrix{2 & 0\cr
                                          0 &  -1\cr}\right),\right.   \lb{4.76}\\
&& \frac{\pt}{\pt e}K_{\mu,e}(t)\left|_{(\mu,e)=(\mu,0)}
       = {-\cos t}\left(\matrix{1+2\mu^* & 0 \cr
                                        0 & 1-\mu^*\cr}\right).\right.   \lb{4.77}\eea
Therefore from (\ref{4.71}) and (\ref{4.74})-(\ref{4.77}) we have
\bea  \<\frac{\pt}{\pt\mu}A(\mu^*,0)x_0,x_0\>
&=& \<\frac{\pt}{\pt\mu}K_{\mu^*,0}z,z\>    \nn\\
&=& \int_0^{2\pi}[2\tilde{a}_0^2\sin^2\frac{t}{2}
                -\cos^2\frac{t}{2}]dt  \nn\\
&=& \pi(2\tilde{a}_0^2-1),  \nn\\
&=& \pi\frac{97-15\sqrt{97}}{64}   \lb{4.78}\eea
and
\bea  \<\frac{\pt}{\pt e}A(\mu^*,0)x_0,x_0\>
&=& \<\frac{\pt}{\pt e}K_{\mu^*,0}z,z\>   \nn\\
&=& -\int_0^{2\pi}[(1+2\mu^*)\tilde{a}_0^2\cos{t}\sin^2\frac{t}{2}
                +(1-\mu^*)\cos{t}\cos^2\frac{t}{2}]dt  \nn\\
&=& \pi\frac{(1+2\mu^*)\tilde{a}_0^2-(1-\mu^*)}{2}.  \nn\\
&=&\pi\frac{-33+15\sqrt{97}}{1024}.    \lb{4.79}\eea
Therefore by (\ref{4.73}) and (\ref{4.78})-(\ref{4.79}),
we obtain
\be  \mu'(0) = \frac{291+15\sqrt{97}}{3104}.  \lb{4.80}\ee
The other tangent can be compute similarly.
Thus the theorem is proved.\hb

Now we can give the

{\bf Proof of the first half of Theorem 1.3. }
Here we give proofs for items (i)-(iii) of this theorem.

{\bf (i)} By Theorem \ref{Th:om.degenerate.curves}, for $\om=-1$,
$\mu_i(e,-1)$ is real analytic on $e\in[0,1)$ for $i=1,2$.

{\bf (ii)} By the computations in Section 3.2,
the only $-1$-degenerate point in the segment $[0,1]\times\{0\}$ is $(\mu,e)=(\mu^*,0)$,
which is a two-fold $-1$-degenerate point.
Because these two curves bifurcate out from $(\mu^*,0)$ in different angles
with tangents $\pm\frac{291+15\sqrt{97}}{3104}$ respectively when $e>0$ is small
by Theorem \ref{Th:tangent.direction},
they are different from each other at least near $(\mu^*,0)$.
Because of analyticity, the intersection points of these two curves can only be isolated.
That $\lim{\mu_i(e,-1)}\to1$ as $e\to1$ for $i=1,2$
follows by the similar arguments in the Section 5 of \cite{HLS}.

{\bf (iii)}  It follows from the computations in Section 3.2,
Lemma \ref{Lemma:increasing.of.index} and Theorem \ref{Th:om.degenerate.curves}.
\hb

\subsection{The hyperbolic region and the symplectic normal forms of $\ga_{\mu,e}(2\pi)$}

For every $e\in[0,1)$, we recall
$$
\mu_l(e)=\inf\{\mu'\in[0,1]|\sg(\ga_{\mu,e}(2\pi))\cap\U\ne\emptyset,\;\forall \mu\in[0,\mu']\},
$$
and
$$
\Ga_l=\{(\mu_l(e),e)\in[0,1]\times[0,1)|e\in[0,1)\}.
$$
By similar arguments of Lemma 9.1 and Corollary 9.2 in \cite{HLS}, we have
\begin{lemma}\lb{Lm:hyperbolic.region}
(i) If $0\le\mu_1<\mu_2\le1$ and $\ga_{\mu_2,e}(2\pi)$ is hyperbolic, so does $\ga_{\mu_1,e}(2\pi)$.
Consequently, the hyperbolic region of $\ga_{\mu,e}(2\pi)$ in $[0,1]\times[0,1)$ is connected.

(ii) For any fixed $e\in[0,1)$, every matrix $\ga_{\mu,e}(2\pi)$ is hyperbolic if $0<\mu<\mu_l(e)$
for $\mu_l(e)$ defined by (\ref{mu_l}).

(iii) We have
\be
\sum_{\mu\in[\mu_l(e),1]}\nu_\om(\ga_{\mu,e}(2\pi))=2,\quad \forall\om\in\U\backslash\{1\}.
\ee

(iv) For every $e\in[0,1)$, we have
\be
\sum_{\mu\in[0,\mu_m(e))}\nu_{-1}(\ga_{\mu,e}(2\pi))=0,\quad
\sum_{\mu\in[\mu_m(e),1]}\nu_{-1}(\ga_{\mu,e}(2\pi))=2.
\ee
\end{lemma}

Now we can give

{\bf The Proof of the second half of Theorem 1.3. }Here we give proofs for items (iv)-(x) of this theorem.
Some arguments below are use the methods in the proof of Theorem 1.2 in \cite{HLS}.

{\bf (iv)} It follows from Lemma \ref{Lm:hyperbolic.region}(ii).

{\bf (v)} In fact, if the function $\mu_l(e)$ is not continuous in $e\in[0,1)$,
then there exist some $\hat{e}\in[0,1)$,
a sequence $\{e_i|i\in\N\}\backslash\{\hat{e}\}$ and $\mu_0\in[0,1]$ such that
\be\lb{non.continuous.assumption}
\mu_l(e_i)\rightarrow\mu_0\ne\mu_l(\hat{e})\quad and\quad e_i\rightarrow\hat{e}\quad
as\quad i\rightarrow+\infty.
\ee
We continue in two cases according to the sign of the difference $\mu_0-\mu_l(\hat{e})$.

On the one hand, by the definition of $\mu_l(e_i)$ we have $\sigma(\ga_{\mu_l(e_i),e_i}(2\pi))\cap\U\ne\emptyset$
for every $e_i$. By the continuity of eigenvalues of $\ga_{\mu_l(e_i),e_i}(2\pi)$ in $i$
 and (\ref{non.continuous.assumption}),
we obtain
\be
\sigma(\ga_{\mu_0,\hat{e}}(2\pi))\cap\U\ne\emptyset.
\nn
\ee
Thus by Lemma \ref{Lm:hyperbolic.region}, this would yield a contradiction if $\mu_0<\mu_l(\hat{e})$.

On the other hand, we suppose $\mu_0>\mu_l(\hat{e})$.
By Lemma \ref{Lm:hyperbolic.region}, for all $i\ge1$, we have
\be\lb{hyperbolic}
\sigma(\ga_{\mu,e_i}(2\pi))\cap\U=\emptyset,\quad \forall\mu\in(0,\mu_l(e_i)).
\ee
Then by the continuity of $\mu_m(e)$ in $e$, (\ref{non.continuous.assumption}) and (\ref{hyperbolic}),
we obtain
\be
\mu_l(\hat{e})<\mu_0\le\mu_m(\hat{e}).
\ee
Let $\om_0\in\sigma(\ga_{\mu_l(\hat{e}),\hat{e}}(2\pi))\cap\U$, which exists by the definition of $\mu_l(\hat{e})$.

Moreover, let $L=\{(\mu,\hat{e})|\mu\in(0,\mu_l(\hat{e}))\}$
and $L_i=\{(\mu,e_i)|\mu\in(0,\mu_l(e_i))\}$ for all $i\ge1$.
Note that by (\ref{index.of.0e}), Lemma \ref{Lemma:increasing.of.index}(iii) and Lemma \ref{Lm:hyperbolic.region}, we obtain
\be\lb{4.29}
i_{\om_0}(\ga_{\mu,e})=\nu_{\om_0}(\ga_{\mu,e})=0,\quad
\forall(\mu,e)\in L\cup \bigcup_{i\ge1}{L_i}.
\ee
Specially, we have
\be
i_{\om_0}(\ga_{\mu_l(\hat{e}),\hat{e}})=0,\quad
\nu_{\om_0}(\ga_{\mu_l(\hat{e}),\hat{e}})\ge1.
\ee
Therefore by Lemma \ref{Lemma:increasing.of.index}(iii) and the definition of $\om_0$,
there exist $\hat\mu\in(\mu_l(\hat{e}),\mu_0)$ sufficiently close to $\mu_l(\hat{e})$ such that
\be\lb{4.31}
i_{\om_0}(\ga_{\hat\mu,\hat{e}})= i_{\om_0}(\ga_{\mu_l(\hat{e}),\hat{e}})
+\nu_{\om_0}(\ga_{\mu_l(\hat{e}),\hat{e}})\ge1.
\ee
This estimate (\ref{4.31}) in facts holds for all $\mu\in(\mu_l(\hat{e}),\hat\mu]$ too.
Note that $(\hat\mu,\hat{e})$ is an accumulation point of $\cup_{i\ge1}{L_i}$.
Consequently for each $i\ge1$, there exist $(\mu_i,e_i)\in L_i$ such that
$\ga_{\mu_i,e_i}\in\P_{2\pi}(4)$ is $\om_0$ non-degenerate, $\mu_i\to\hat\mu$ in $\R$,
and $\ga_{\mu_i,e_i}\to\ga_{\hat\mu,\hat e}$ in $\P_{2\pi}(4)$ as $i\to\infty$.
Therefore by (\ref{4.29}), (\ref{4.31}), the Definition 5.4.2 of the $\om_0$-index of
$\om_0$-degenerate path $\ga_{\hat\mu,\hat e}$ on p.129 and Theorem6.1.8 on p.142 of \cite{Lon4},
we obtain the following contradiction
\be
1\le i_{\om_0}(\ga_{\hat\mu,\hat e})\le i_{\om_0}(\ga_{\mu_i,e_i})=0,
\ee
for $i\ge1$ large enough. Thus the continuity of $\mu_l(e)$ in $e\in[0,1)$ is proved.

Now we prove the claim $\lim_{e\to1}\mu_l(e)=1$.
We argue by contradiction, and suppose there exist $e_i\to1$ as $i\to+\infty$ such that
$\lim_{e\to1}\mu_l(e)=\mu_0$ for some $0\le\mu_0<1$.
Then at least one of the following two cases must occur:
{\bf (A)} There exists a subsequence $\hat{e}_i$ of $e_i$ such that $\mu_l(\hat{e}_{i+1})\le\mu_l(\hat{e}_{i})$
for all $i\in\N$;
{\bf (B)} There exists a subsequence $\hat{e}_i$ of $e_i$ such that $\mu_l(\hat{e}_{i+1})\ge\mu_l(\hat{e}_{i})$
for all $i\in\N$.

If Case (A) happens, for this $\mu_0$, by a similar argument of Theorem 1.7 in \cite{HLS},
there exists $e_0>0$ sufficiently close to $1$ such that
$\ga_{\mu,e}(2\pi)$ is hyperbolic for all $(\mu,e)$ in the region $(0,\mu_l(\hat{e}_i)]\times[e_0,1)$.
Then by the monotonicity of Case (A) we obtain
\be
\mu_0\le\mu_l(\hat{e}_{i+m})\le\mu_l(\hat{e}_{i}),\quad
\forall m\in\N.
\ee
Therefore $(\mu_l(\hat{e}_{i+m}),\hat{e}_{i+m})$ will get into this region for sufficiently large $m\ge1$,
which contract to the definition of $\mu_l(\hat{e}_{i+m})$.

If Case (B) happens, the proof is similar. Thus (v) holds.

{\bf (vi)}  By our study in Section 3.2, we have $({8\over9},0)\in\Ga_l\backslash\Ga_m$.
Thus there exist an $\tilde{e}\in(0,1]$ such that $\mu_l(e)<\mu_m(e)$ for all $e\in[0,\tilde{e})$.
Therefore, $\Ga_l$ and $\Ga_m$ are different curves.

{\bf (vii)} If $\mu_l(e)<\mu<\mu_m(e)$, then by the definitions of the degenerate curves and
Lemma \ref{Lemma:increasing.of.index} $(iii)$, we have
\begin{equation}
i_1(\ga_{\mu,e})=0,\quad\quad \nu_1(\ga_{\mu,e})=0, \lb{4.89}
\end{equation}
and
\begin{equation}
i_{-1}(\ga_{\mu,e})=0,\quad\quad \nu_{-1}(\ga_{\mu,e})=0.  \lb{4.90}
\end{equation}

Assume $\ga_{\mu,e}(2\pi)\approx N_2(e^{\sqrt{-1}\theta},b)$ for some $\theta\in(0,\pi)\cup(\pi,2\pi)$.
Without lose of generality, we suppose $\theta\in(0,\pi)$.
Let $\om_0=e^{\sqrt{-1}\theta}$, we have $\nu_{\om_0}(\ga_{\mu,e}(2\pi))\ge1$.
Then for any $\om\in\U,\om\ne\om_0$, we have
\be
i_{\om}(\ga_{\mu,e})=i_1(\ga_{\mu,e})=0
\ee
or
\be
i_{\om}(\ga_{\mu,e})=i_1(\ga_{\mu,e})-S_{N_2(e^{\sqrt{-1}\theta},b)}^-(e^{\sqrt{-1}\theta})+S_{N_2(e^{\sqrt{-1}\theta},b)}^+(e^{\sqrt{-1}\theta})=0.
\ee
Then by the sub-continuous of $i_\om(\ga_{\mu,e})$ with respect to $\om$,
we have $i_{\om}(\ga_{\mu,e})=0,\; \forall\om\in\U$.
Moreover, by Corollary \ref{C4.5}, we have
\be
i_{\om}(\ga_{\tilde\mu,e})=0,\qquad \forall\om\in\U, \tilde\mu\in(0,\mu).
\ee
Therefore, by the definition of $\mu_l(e)$ of (\ref{mu_l}), we have $\mu_l(e)\ge\mu$.
It contradicts $\mu_l(e)<\mu<\mu_m(e)$.

Then we can suppose $\ga_{\mu,e}(2\pi)\approx M_1\diamond M_2$
where $M_1$ and $M_2$ are two basic normal forms in $\Sp(2)$
defined in Section 5.2 below. Let $\ga_1$ and $\ga_2$ be two paths in $\P_{2\pi}(2)$ such that $\ga_1(2\pi)=M_1$,
$\ga_2(2\pi)=M_2$ and $\ga_{\mu,e}\sim\ga_1\diamond\ga_2$. Then
\begin{equation}
0=i_1(\ga_{\mu,e})=i_1(\ga_1)+i_1(\ga_2).
\end{equation}
By the definition of $\mu_k(s)$, $M_1$ and $M_2$ cannot be both hyperbolic,
and without loose of generality, we suppose $M_1=R(\theta_1)$.
Then $i_1(\ga_1)$ is odd, and hence $i_1(\ga_2)$ is also odd.
By Theorem 4 to Theorem 7 of Chapter 8 on pp.179-183
in \cite{Lon4} and using notations there,
we must have $M_2=D(-2)$ or $M_2=R(\theta_2)$ for some $\theta_2\in(0,\pi)\cup(\pi,2\pi)$.

If $M_2=D(-2)$, then we have $i_{-1}(\ga_1)-i_1(\ga_1)=\pm1$ and $i_{-1}(\ga_2)-i_1(\ga_2)=0$.
Therefore $i_{-1}(\ga_{\mu,e}(2\pi))=i_{-1}(\ga_1)+i_{-1}(\ga_2)$
and $i_{1}(\ga_{\mu,e}(2\pi))=i_{1}(\ga_1)+i_{1}(\ga_2)$ has the different odevity,
which contradicts (\ref{4.89}) and (\ref{4.90}).
Then we have $M_2=R(\theta_2)$.

Moreover, if $\th_1\in(\pi,2\pi)$, we must have $\th_2\in(0,\pi)$,
otherwise $i_{-1}(\ga_1)-i_1(\ga_1)=1$ and $i_{-1}(\ga_2)-i_1(\ga_2)=1$ and hence
\be
i_{-1}(\ga_{\mu,e})=i_{-1}(\ga_1)+i_{-1}(\ga_2)
=i_1(\ga_1)+i_1(\ga_1)+2=2,
\ee
which contradicts (\ref{4.90}).
Similarly, if if $\th_1\in(0,\pi)$, we must have $\th_2\in(\pi,2\pi)$.

{\bf (viii)}   If $\mu_m(e)<\mu<\mu_r(e)$, then by the definitions of the degenerate curves and
Lemma \ref{Lemma:increasing.of.index} $(iii)$, we have
\begin{equation}
i_1(\ga_{\mu,e})=0,\quad\quad \nu_1(\ga_{\mu,e})=0, \lb{1.index}
\end{equation}
and
\begin{equation}
i_{-1}(\ga_{\mu,e})=1,\quad\quad \nu_{-1}(\ga_{\mu,e})=0.  \lb{-1.index}
\end{equation}

Firstly, if $\ga_{\mu,e}(2\pi)\approx N_2(e^{\sqrt{-1}\theta},b)$ for some $\theta\in(0,\pi)\cup(\pi,2\pi)$,
we have
\be
i_{-1}(\ga_{\mu,e})=i_1(\ga_{\bb,e})-S_{N_2(e^{\sqrt{-1}\theta},b)}^-(e^{\sqrt{-1}\theta})+S_{N_2(e^{\sqrt{-1}\theta},b)}^+(e^{\sqrt{-1}\theta})=i_1(\ga_{\mu,e})
\ee
or
\be
i_{-1}(\ga_{\mu,e})=i_1(\ga_{\mu,e})-S_{N_2(e^{\sqrt{-1}\theta},b)}^-(e^{\sqrt{-1}(2\pi-\theta)})+S_{N_2(e^{\sqrt{-1}\theta},b)}^+(e^{\sqrt{-1}(2\pi-\theta)})=i_1(\ga_{\mu,e}),
\ee
which contradicts (\ref{1.index}) and (\ref{-1.index}).

Then we can suppose $\ga_{\mu,e}(2\pi)\approx M_1\diamond M_2$ where $M_1$ and $M_2$ are two basic normal forms in $\Sp(2)$
defined in Section 5.2 below. Let $\ga_1$ and $\ga_2$ be two paths in $\P_{2\pi}(2)$ such that $\ga_1(2\pi)=M_1$,
$\ga_2(2\pi)=M_2$ and $\ga_{\mu,e}\sim\ga_1\diamond\ga_2$. Then
\begin{equation}
0=i_1(\ga_{\mu,e})=i_1(\ga_1)+i_1(\ga_2).
\end{equation}
By the definition of $\mu_k(s)$, $M_1$ and $M_2$ cannot be both hyperbolic,
and without loose of generality, we suppose $M_1=R(\theta_1)$.
Then $i_1(\ga_1)$ is odd, and hence $i_1(\ga_2)$ is also odd.
By Theorem 4 to Theorem 7 of Chapter 8 on pp.179-183
in \cite{Lon4} and using notations there,
we must have $M_2=D(-2)$ or $M_2=R(\theta_2)$ for some $\theta_2\in(0,\pi)\cup(\pi,2\pi)$.

If $M_2=R(\theta_2)$, then we have $i_{-1}(\ga_1)-i_1(\ga_1)=\pm1$ and $i_{-1}(\ga_2)-i_1(\ga_2)=\pm1$.
Therefore $i_{-1}(\ga_{\mu,e}(2\pi))=i_{-1}(\ga_1)+i_{-1}(\ga_2)$
and $i_{1}(\ga_{\mu,e}(2\pi))=i_{1}(\ga_1)+i_{1}(\ga_2)$ has the same odevity,
which contradicts to (\ref{1.index}) and (\ref{-1.index}).
Then we have $M_2=D(-2)$.

Moreover, if $\th_1\in(0,\pi)$, we have
\be
i_{-1}(\ga_{\mu,e})=i_{-1}(\ga_1)+i_{-1}(\ga_2)=
i_1(\ga_1)-S_{R(\theta_1)}^-(e^{\sqrt{-1}\theta})+S_{R(\theta_2)}^+(e^{\sqrt{-1}\theta})+i_{-1}(\ga_2)
=i_1(\ga_1)-1+i_1(\ga_1)=-1,
\ee
which contradicts (\ref{-1.index}).
Thus (viii) is proved. (ix) can be proved by similar steps.
\hb

\begin{remark}
Remark \ref{R1.4} can be obtained by similar arguments to Theorem 1.3 (vii)-(ix).
\end{remark}

\setcounter{equation}{0}
\section{Appendix: $\omega$-Maslov-type indices and $\omega$-Morse indices}

Let $(\R^{2n},\Omega)$ be the standard symplectic vector space with coordinates
$(x_1,...,x_n,y_1,...,y_n)$ and the symplectic form $\Omega=\sum_{i=1}^{n}dx_i \wedge dy_i$.
Let $J=\left(\matrix{0&-I_n\cr
                 I_n&0\cr}\right)$ be the standard symplectic matrix, where $I_n$
is the identity matrix on $\R^n$.

As usual, the symplectic group $\Sp(2n)$ is defined by
$$ \Sp(2n) = \{M\in {\rm GL}(2n,\R)\,|\,M^TJM=J\}, $$
whose topology is induced from that of $\R^{4n^2}$. For $\tau>0$
we are interested in paths in $\Sp(2n)$:
$$ \P_{\tau}(2n) = \{\ga\in C([0,\tau],\Sp(2n))\,|\,\ga(0)=I_{2n}\}, $$
which is equipped with the topology induced from that of $\Sp(2n)$.
For any $\om\in\U$ and $M\in\Sp(2n)$, the following real function was
introduced in \cite{Lon2}:
$$ D_{\om}(M) = (-1)^{n-1}\ol{\om}^n\det(M-\om I_{2n}). $$
Thus for any $\om\in\U$ the following codimension $1$ hypersurface
in $\Sp(2n)$ is defined (\cite{Lon2}):
$$ \Sp(2n)_{\om}^0 = \{M\in\Sp(2n)\,|\, D_{\om}(M)=0\}.  $$
For any $M\in \Sp(2n)_{\om}^0$, we define a co-orientation of
$\Sp(2n)_{\om}^0$ at $M$ by the positive direction
$\frac{d}{dt}Me^{t J}|_{t=0}$ of the path $Me^{t J}$ with $0\le t\le
\varepsilon$ and $\varepsilon$ being a small enough positive number. Let
\bea
\Sp(2n)_{\om}^{\ast} &=& \Sp(2n)\bs \Sp(2n)_{\om}^0,   \nn\\
\P_{\tau,\om}^{\ast}(2n) &=&
      \{\ga\in\P_{\tau}(2n)\,|\,\ga(\tau)\in\Sp(2n)_{\om}^{\ast}\}, \nn\\
\P_{\tau,\om}^0(2n) &=& \P_{\tau}(2n)\bs \P_{\tau,\om}^{\ast}(2n). \nn\eea
For any two continuous paths $\xi$ and $\eta:[0,\tau]\to\Sp(2n)$ with
$\xi(\tau)=\eta(0)$, we define their concatenation by:
$$ \eta\ast\xi(t) = \left\{\matrix{
            \xi(2t), & \quad {\rm if}\;0\le t\le \tau/2, \cr
            \eta(2t-\tau), & \quad {\rm if}\; \tau/2\le t\le \tau. \cr}\right. $$
Given any two $2m_k\times 2m_k$ matrices of square block form
$M_k=\left(\matrix{A_k&B_k\cr
                                C_k&D_k\cr}\right)$ with $k=1, 2$,
the symplectic sum of $M_1$ and $M_2$ is defined (cf. \cite{Lon2} and \cite{Lon4}) by
the following $2(m_1+m_2)\times 2(m_1+m_2)$ matrix $M_1\dm M_2$:
\be M_1\dm M_2=\left(\matrix{A_1 &   0 & B_1 &   0\cr
                             0   & A_2 &   0 & B_2\cr
                             C_1 &   0 & D_1 &   0\cr
                             0   & C_2 &   0 & D_2\cr}\right),   \lb{2.1} \ee
and $M^{\dm k}$ denotes the $k$ copy $\dm$-sum of $M$. For any two paths $\ga_j\in\P_{\tau}(2n_j)$
with $j=0$ and $1$, let $\ga_0\dm\ga_1(t)= \ga_0(t)\dm\ga_1(t)$ for all $t\in [0,\tau]$.

As in \cite{Lon4}, for $\lm\in\R\bs\{0\}$, $a\in\R$, $\th\in (0,\pi)\cup (\pi,2\pi)$,
$b=\left(\matrix{b_1 & b_2\cr
                 b_3 & b_4\cr}\right)$ with $b_i\in\R$ for $i=1, \ldots, 4$, and $c_j\in\R$
for $j=1, 2$, we denote respectively some normal forms by
\bea
&& D(\lm)=\left(\matrix{\lm & 0\cr
                         0  & \lm^{-1}\cr}\right), \qquad
   R(\th)=\left(\matrix{\cos\th & -\sin\th\cr
                        \sin\th  & \cos\th\cr}\right),  \nn\\
&& N_1(\lm, a)=\left(\matrix{\lm & a\cr
                             0   & \lm\cr}\right), \qquad
   N_2(e^{\sqrt{-1}\th},b) = \left(\matrix{R(\th) & b\cr
                                           0      & R(\th)\cr}\right),  \nn\\
&& M_2(\lm,c)=\left(\matrix{\lm &   1 &       c_1 &         0 \cr
                              0 & \lm &       c_2 & (-\lm)c_2 \cr
                              0 &   0 &  \lm^{-1} &         0 \cr
                              0 &   0 & -\lm^{-2} &  \lm^{-1} \cr}\right). \nn\eea
Here $N_2(e^{\sqrt{-1}\th},b)$ is {\bf trivial} if $(b_2-b_3)\sin\th>0$, or {\bf non-trivial}
if $(b_2-b_3)\sin\th<0$, in the sense of Definition 1.8.11 on p.41 of \cite{Lon4}. Note that
by Theorem 1.5.1 on pp.24-25 and (1.4.7)-(1.4.8) on p.18 of \cite{Lon4}, when $\lm=-1$ there hold
\bea
c_2 \not= 0 &{\rm if\;and\;only\;if}\;& \dim\ker(M_2(-1,c)+I)=1, \nn\\
c_2 = 0 &{\rm if\;and\;only\;if}\;& \dim\ker(M_2(-1,c)+I)=2. \nn\eea
Note that we have $N_1(\lm,a)\approx N_1(\lm, a/|a|)$ for $a\in\R\bs\{0\}$ by symplectic coordinate
change, because
$$ \left(\matrix{1/\sqrt{|a|} & 0\cr
                           0  & \sqrt{|a|}\cr}\right)
   \left(\matrix{\lm & a\cr
                  0  & \lm\cr}\right)
   \left(\matrix{\sqrt{|a|} & 0\cr
                           0  & 1/\sqrt{|a|}\cr}\right) = \left(\matrix{\lm & a/|a|\cr
                                                                         0  & \lm\cr}\right). $$

\begin{definition}\lb{D2.1} (\cite{Lon2}, \cite{Lon4})
For any $\om\in\U$ and $M\in \Sp(2n)$, define
\be \nu_{\om}(M)=\dim_{\C}\ker_{\C}(M - \om I_{2n}).  \lb{2.2}\ee

For every $M\in \Sp(2n)$ and $\om\in\U$, as in Definition 1.8.5 on p.38 of \cite{Lon4}, we define the
{\bf $\om$-homotopy set} $\Om_{\om}(M)$ of $M$ in $\Sp(2n)$ by
$$  \Om_{\om}(M)=\{N\in\Sp(2n)\,|\, \nu_{\om}(N)=\nu_{\om}(M)\},  $$
and the {\bf homotopy set} $\Om(M)$ of $M$ in $\Sp(2n)$ by
\bea  \Om(M)=\{N\in\Sp(2n)\,&|&\,\sg(N)\cap\U=\sg(M)\cap\U,\,{\it and}\; \nn\\
         &&\qquad \nu_{\lm}(N)=\nu_{\lm}(M)\qquad\forall\,\lm\in\sg(M)\cap\U\}.  \nn\eea
We denote by $\Om^0(M)$ (or $\Om^0_{\om}(M)$) the path connected component of $\Om(M)$ ($\Om_{\om}(M)$)
which contains $M$, and call it the {\bf homotopy component} (or $\om$-{\bf homtopy component}) of $M$ in
$\Sp(2n)$. Following Definition 5.0.1 on p.111 of \cite{Lon4}, for $\om\in \U$ and $\ga_i\in \P_{\tau}(2n)$
with $i=0, 1$, we write $\ga_0\sim_{\om}\ga_1$ if $\ga_0$ is homotopic to $\ga_1$ via
a homotopy map $h\in C([0,1]\times [0,\tau], \Sp(2n))$ such that $h(0)=\ga_0$, $h(1)=\ga_1$, $h(s)(0)=I$,
and $h(s)(\tau)\in \Om_{\om}^0(\ga_0(\tau))$ for all $s\in [0,1]$. We write also $\ga_0\sim \ga_1$, if
$h(s)(\tau)\in \Om^0(\ga_0(\tau))$ for all $s\in [0,1]$ is further satisfied.
\end{definition}

Following Definition 1.8.9 on p.41 of \cite{Lon4}, we call the above matrices $D(\lm)$, $R(\th)$, $N_1(\lm,a)$
and $N_2(\om,b)$ basic normal forms of symplectic matrices. As proved in \cite{Lon2} and \cite{Lon3} (cf.
Theorem 1.9.3 on p.46 of \cite{Lon4}), every $M\in\Sp(2n)$ has its basic normal form decomposition in $\Om^0(M)$
as a $\dm$-sum of these basic normal forms. This is very important when we derive basic normal forms
for $\ga_{\bb,e}(2\pi)$ to compute the $\om$-index $i_{\om}(\ga_{\bb,e})$ of the path $\ga_{\bb,e}$ later
in this paper.

We define a special continuous symplectic path $\xi_n\subset {\Sp}(2n)$ by
\be \xi_n(t) = \left(\matrix{2-\frac{t}{\tau} & 0 \cr
                             0 &  (2-\frac{t}{\tau})^{-1}\cr}
               \right)^{\dm n} \qquad {\rm for}\;0\le t\le \tau.  \lb{2.3}\ee

\begin{definition} (\cite{Lon2}, \cite{Lon4})\lb{D2.2}
{For any $\tau>0$ and $\ga\in \P_{\tau}(2n)$, define
\be \nu_{\om}(\ga)= \nu_{\om}(\ga(\tau)).  \lb{2.4}\ee

If $\ga\in\P_{\tau,\om}^{\ast}(2n)$, define
\be i_{\om}(\ga) = [\Sp(2n)_{\om}^0: \ga\ast\xi_n],  \lb{2.5}\ee
where the right hand side of (\ref{2.5}) is the usual homotopy intersection number, and
the orientation of $\ga\ast\xi_n$ is its positive time direction under homotopy with
fixed end points.

If $\ga\in\P_{\tau,\om}^0(2n)$, we let $\mathcal{F}(\ga)$ be the set of all open
neighborhoods of $\ga$ in $\P_{\tau}(2n)$, and define
\be i_{\om}(\ga) = \sup_{U\in\mathcal{F}(\ga)}\inf\{i_{\om}(\beta)\,|\,
                       \beta\in U\cap\P_{\tau,\om}^{\ast}(2n)\}.      \lb{2.6}\ee
Then
$$ (i_{\om}(\ga), \nu_{\om}(\ga)) \in \Z\times \{0,1,\ldots,2n\}, $$
is called the index function of $\ga$ at $\om$. }
\end{definition}

\begin{definition} (\cite{Lon2}, \cite{Lon4})\lb{D2.3}
For any $M\in\Sp(2n)$ and $\om\in\U$, choosing $\tau>0$ and $\ga\in\P_\tau(2n)$ with $\ga(\tau)=M$,
we define
\be
S_M^{\pm}(\om)=\lim_{\epsilon\rightarrow0^+}\;i_{\exp(\pm\epsilon\sqrt{-1}\om)}(\ga)-i_\om(\ga).
\ee
They are called the splitting numbers of $M$ at $\om$.
\end{definition}

We refer to \cite{Lon4} for more details on this index theory of symplectic matrix paths
and periodic solutions of Hamiltonian system.

For $T>0$, suppose $x$ is a critical point of the functional
$$ F(x)=\int_0^TL(t,x,\dot{x})dt,  \qquad \forall\,\, x\in W^{1,2}(\R/T\Z,\R^n), $$
where $L\in C^2((\R/T\Z)\times \R^{2n},\R)$ and satisfies the
Legendrian convexity condition $L_{p,p}(t,x,p)>0$. It is well known
that $x$ satisfies the corresponding Euler-Lagrangian
equation:
\bea
&& \frac{d}{dt}L_p(t,x,\dot{x})-L_x(t,x,\dot{x})=0,    \label{2.7}\\
&& x(0)=x(T),  \qquad \dot{x}(0)=\dot{x}(T).    \label{2.8}\eea

For such an extremal loop, define
\bea
P(t) &=& L_{p,p}(t,x(t),\dot{x}(t)),  \nn\\
Q(t) &=& L_{x,p}(t,x(t),\dot{x}(t)),  \nn\\
R(t) &=& L_{x,x}(t,x(t),\dot{x}(t)).  \nn\eea
Note that
\be F\,''(x)=-\frac{d}{dt}(P\frac{d}{dt}+Q)+Q^T\frac{d}{dt}+R. \lb{2.9}\ee

For $\omega\in\U$, set
\be  D(\omega,T)=\{y\in W^{1,2}([0,T],\C^n)\,|\, y(T)=\omega y(0) \}.   \lb{2.10}\ee
We define the $\omega$-Morse index $\phi_\omega(x)$ of $x$ to be the dimension of the
largest negative definite subspace of
$$ \langle F\,''(x)y_1,y_2 \rangle, \qquad \forall\;y_1,y_2\in D(\omega,T), $$
where $\langle\cdot,\cdot\rangle$ is the inner product in $L^2$. For $\omega\in\U$, we
also set
\be  \ol{D}(\omega,T)= \{y\in W^{2,2}([0,T],\C^n)\,|\, y(T)=\omega y(0), \dot{y}(T)=\om\dot{y}(0) \}.
                     \lb{2.11}\ee
Then $F''(x)$ is a self-adjoint operator on $L^2([0,T],\R^n)$ with domain $\ol{D}(\omega,T)$.
We also define
\[\nu_\omega(x)=\dim\ker(F''(x)).\]

In general, for a self-adjoint operator $A$ on the Hilbert space $\mathscr{H}$, we set
$\nu(A)=\dim\ker(A)$ and denote by $\phi(A)$ its Morse index which is the maximum dimension
of the negative definite subspace of the symmetric form $\langle A\cdot,\cdot\rangle$. Note
that the Morse index of $A$  is equal to the total multiplicity of the negative eigenvalues
of $A$.

On the other hand, $\td{x}(t)=(\partial L/\partial\dot{x}(t),x(t))^T$ is the solution of the
corresponding Hamiltonian system of (\ref{2.7})-(\ref{2.8}), and its fundamental solution
$\gamma(t)$ is given by
\bea \dot{\gamma}(t) &=& JB(t)\gamma(t),  \lb{2.12}\\
     \gamma(0) &=& I_{2n},  \lb{2.13}\eea
with
\be B(t)=\left(\matrix{P^{-1}(t)& -P^{-1}(t)Q(t)\cr
                       -Q(t)^TP^{-1}(t)& Q(t)^TP^{-1}(t)Q(t)-R(t)\cr}\right). \lb{2.14}\ee

\begin{lemma}(Long, \cite{Lon4}, p.172)\lb{L2.3}  
For the $\omega$-Morse index $\phi_\omega(x)$ and nullity $\nu_\omega(x)$ of the solution $x=x(t)$
and the $\omega$-Maslov-type index $i_\omega(\gamma)$ and nullity $\nu_\omega(\gamma)$ of the symplectic
path $\ga$ corresponding to $\td{x}$, for any $\omega\in\U$ we have
\be \phi_\omega(x) = i_\omega(\gamma), \qquad \nu_\omega(x) = \nu_\omega(\gamma).  \lb{2.15}\ee
\end{lemma}

A generalization of the above lemma to arbitrary  boundary conditions is given in \cite{HS1}.
For more information on these topics, we refer to \cite{Lon4}.

\medskip

\noindent {\bf Acknowledgements.}
The authors thank sincerely Professor Yiming Long for
his precious help and useful discussions.

\end{document}